\documentstyle[12pt]{article}
\newcommand{\ra}{\rightarrow}
\newtheorem{theorem}{Theorem}
\newtheorem{lemma}[theorem]{Lemma}
\newtheorem{proposition}[theorem]{Proposition}

\newtheorem{corollary}[theorem]{Corollary}
\newtheorem{remar}[theorem]{Remark}
\newenvironment{proof}{Proof:\ \ \ }{\QED}
\newenvironment{remark}{\begin{remar}\rm}{\end{remar}}
\newcommand{\lex}{\,\vec{\amalg}\,}
\newcommand{\QED}{{\unskip\nobreak\hfil\penalty50%
\hskip1em\hbox{}\nobreak\hfil $\Box$%
\parfillskip=0pt \finalhyphendemerits=0 \par\medskip\noindent}}

\newcommand{\n}{\par\noindent}
\newcommand{\card}{\mbox{\rm card}\,}
\newcommand{\Mon}{\mbox{\rm Mon}\,}
\newcommand{\Aut}{\mbox{\rm Aut}\,}
\newcommand{\eqad}{\stackrel{+}{\sim}}
\newcommand{\Llra}{\Longleftrightarrow}
\newcommand{\eqmult}{\stackrel{\cdot}{\sim}}
\newcommand{\Neg}{\mbox{\rm Neg}\,}
\newcommand{\support}{\mbox{\rm support}\,}
\newcommand{\sn}{\par\smallskip\noindent}
\newcommand{\mn}{\par\medskip\noindent}
\newcommand{\bn}{\par\bigskip\noindent}

\font\tenlv=msbm10 scaled 1200
\font\sevenlv=msbm7 scaled 1200
\font\fivelv=msbm5 scaled 1200
\def\lv #1{{\mathchoice{{\hbox{\tenlv #1}}}{{\hbox{\tenlv #1}}}
{{\hbox{\sevenlv #1}}}{{\hbox{\fivelv #1}}}}}
\newcommand{\N}{\lv N}
\newcommand{\Q}{\lv Q}
\newcommand{\R}{\lv R}
\newcommand{\Z}{\lv Z}
\font\tenlv=msbm10 scaled 1200
\font\sevenlv=msbm7 scaled 1200
\font\fivelv=msbm5 scaled 1200
\def\lv #1{{\mathchoice{{\hbox{\tenlv #1}}}{{\hbox{\tenlv #1}}}
{{\hbox{\sevenlv #1}}}{{\hbox{\fivelv #1}}}}}
\renewcommand{\N}{\lv N}
\renewcommand{\Z}{\lv Z}
\renewcommand{\Q}{\lv Q}
\renewcommand{\R}{\lv R}
\begin{document}            
\title{
$\kappa$-bounded Exponential-Logarithmic
                      Power Series Fields
\footnote{
\n First author partially supported
by an NSERC research grant.
This paper was written while the
first author was on sabbatical leave at Universit\'e Paris 7. The
author wishes to thank the Equipe de Logique de Paris 7 for its support
and hospitality.\n The second author would like to
thank the Israel Science Foundation
   for partial support of this research. {\bf Publication 857}.}}
\author{Salma Kuhlmann and Saharon
Shelah\\[.3cm]}
\date{17.\ 04.\ 2005}
\maketitle
%
\begin{abstract}\noindent
In [K--K--S] it was shown that fields of generalized power series
cannot admit an exponential function. In this paper, we construct
fields of generalized power series with {\it bounded support} which
admit an exponential. We give a natural definition of an
exponential, which makes these fields into models of real
exponentiation.
 The method allows to construct for every $\kappa$ regular
uncountable cardinal,
 $2^{\kappa}$ pairwise non-isomorphic models of real exponentiation (of
cardinality
$\kappa$), but all isomorphic as ordered fields. Indeed, the $2^{\kappa}$
exponentials constructed have pairwise distinct {\it growth rates}. This
method relies on constructing lexicographic chains with many
automorphisms.
\end{abstract}
\section{Introduction.}
In [T], Tarski proved his celebrated result that the elementary theory of the ordered field of real numbers admits elimination of quantifiers,
and gave a recursive axiomatization of its class of models (the class of real closed fields). He asked whether analogous results hold for the elementary theory T$_{\exp}$ of $(\R, \exp)$ (the ordered field of real numbers with {\it exponentiation}). Addressing Tarski's problem, Wilkie [W] established that T$_{\exp}$ is model complete and o-minimal.
Due to these results, the problem of
constructing non-archimedean models of T$_{\exp}$ gained much interest.
\sn
Non-archimedean real closed fields are easy to construct; for example, any field of generalized power series (see Section 2)
$\R((G))$ with exponents in a {\it divisible} ordered abelian group $G\not= 0$
is such a model. However,
in [K--K--S] it was shown that fields of generalized power series
cannot admit an exponential function, so different methods were needed to
construct non-archimedean real closed exponential fields.
In [D--M--M2], van den Dries, Macintyre and Marker construct
non-archimedean models (the logarithmic-exponential power series fields) of T$_{\exp}$ with many interesting properties. In
[K], the
exponential-logarithmic power series fields are constructed, providing
yet another class of models. Although the two construction procedures are different (and produce different models, see [K--T]), both logarithmic-exponential or exponential-logarithmic series models are obtained as countable increasing unions of fields of generalized power series.
In both cases, a partial exponential (logarithm) is constructed on every
member of this union, and the exponential on the union is given by an inductive definition.
\sn
In this paper, we describe a different construction, which offers
several advantages. The procedure is straightforward: we start with any non-empty chain $\Gamma _0$. For a given regular uncountable cardinal $\kappa$,
we form the (uniquely determined) $\kappa$-th iterated lexicographic
power $(\Gamma _{\kappa}, \iota_{\kappa})$ of $\Gamma_0$
(see Section 4).
We take $G_{\kappa}$ and $\R((G_{\kappa}))_{\kappa}$ to be the corresponding
$\kappa$-bounded Hahn group and $\kappa$-bounded power series field respectively (see Section 2). The logarithm on the positive elements of $\R((G_{\kappa}))_{\kappa}$ is now defined by a {\it uniform} formula (\ref {logall}). Under the additional hypothesis that $\kappa = \kappa ^{< \kappa}$,
$\R((G_{\kappa}))_{\kappa}$ is a model of cardinality $\kappa$.
\sn
As application, we construct $2^{\kappa}$ pairwise non-isomorphic models of T$_{\exp}$ (of cardinality $\kappa$), but all isomorphic as real closed fields. This answers a question of D.\ Marker, and establishes an exponential analogue to the main result of [A-K].
\sn
The structure of the paper is as follows. In Section 2, we recall some preliminary notions and facts. In Section 3, we state and prove the Main Lemma: it provides sufficient conditions on a chain $\Gamma$, which allow a uniform definition of
a logarithm on $\R((G_{\kappa}))_{\kappa}$. In Section 4, we give a canonical procedure to obtain chains satisfying the conditions of the Main Lemma. In Proposition 4, an additional sufficient condition, which allows to obtain logarithms satisfying the {\it growth axiom scheme} is given. In Section 5, we complete the construction of the model (Theorem 7). In Section 6, we introduce the {\it logarithmic rank}, which is an isomorphism invariant for the logarithm. Theorem 8 relates the logarithmic rank of our
model to the orbital behaviour of automorphisms of our initial chain $\Gamma_0$. In Section 7, we construct chains with many automorphisms, which in turn
allows the construction of models of T$_{\exp}$ with many logarithms (Theorem 9).
\mn
We would like to thank D. Marker for asking us this question, and T. Green for proof-reading preliminary versions of this paper.
\sn
\section{Preliminaries}
We first need some definitions and general facts.
Let $\Gamma$ be a {\bf chain} (that is, a totally ordered set).
Let $X$, $Y$ be subsets of $\Gamma$.
We write $X<Y$ if $x < y$ for all $x \in X$ and
$y \in Y$. A Dedekind cut in $\Gamma$ is a pair $(X,Y)$ of
disjoint
nonempty convex subsets of $\Gamma$ whose union is $\Gamma$ and $X<Y$. A
Dedekind
cut is a {\bf gap} in $\Gamma$ if $X$ has no last element and
$Y$ has no
first element.
$\Gamma$ is said to be
Dedekind complete if there are no gaps
in $\Gamma$.
We denote by
$\overline{\Gamma}$ the
Dedekind completion of a chain $\Gamma$.
We
say that a point $\alpha\in \Gamma$ has {\bf left character} $\aleph_0$
if
$\{\alpha'\in \Gamma\,;\,\alpha'<\alpha\}$ has cofinality $\aleph_0$,
and dually for right character.
Similarly, the characters of a gap
$\overline{s}$ in a chain $\Gamma$ are those of $\overline{s}$
considered as a point in $\overline{\Gamma}$. If both characters are
$\aleph_0$,
we shall call it an $\aleph_0 \aleph_0$-gap.
\sn
Given chains $\Gamma$ and $\Gamma'$, we denote by
$\Gamma\lex
\Gamma'$
the chain obtained by lexicographically ordering the Cartesian product
$\Gamma\times
\Gamma'$. In other words, we obtain the ordered sum of chains
$\Gamma \lex \Gamma'\simeq\sum _{\gamma \in \Gamma} \Gamma'_{\gamma}$
(where $\Gamma'_{\gamma}$ denotes the $\gamma$-th
copy of $\Gamma'$).
\sn
Let $G$ be a totally ordered abelian group. The
 archimedean
equivalence relation on $G$ is
 defined as follows:
\[\mbox{ For } x,y\in G\setminus \{0\} :\mbox{ } x\eqad y\mbox { if }
\exists n\in
\N
\mbox{ s.t. }n|x| \geq |y|
\mbox{ and } n|y|\geq |x|\;\]
where $|x|:=\max\{x,-x\}$. We set $x << y$ if for all $n\in\N$,  $n|x| < |y|$. We denote by $[x]$ is the archimedean equivalence class of $x$. We totally order the set of archimedean classes as follows: $[y] < [x]$ if $x << y$.
 \mn
Let $(K,+,\cdot,0,1,<)$ be an ordered field. Using the
archimedean
equivalence relation on the ordered abelian group
$(K,+,0,<)$, we can endow $K$
with the
{\bf natural
valuation}
 $v$: for $x,y\in K$, $x, y\not= 0$ define $v(x):=[x]$ and $[x] + [y] :=
[xy]$.
We call $v(K):=\{v(x)\>|\> x\in K, x\not= 0\}$ the {\bf value group},
$\;R_v := \{x\mid x\in K \mbox{ and } v(x)\geq 0\}$ the {\bf valuation
ring},
$\;I_v := \{x\mid x\in K \mbox{ and } v(x)> 0\}\;$
 the {\bf valuation ideal}
(the unique maximal ideal of $R_v$),
$\;U_v^{>0}:=\{x\mid x\in R_v, x>0, v(x)=0\}\;$ {\bf the group of
positive units} of $R_v$. The {\bf residue field} is $\overline
{K}:=R_v/I_v$.
 For  $x,y\in K^{>0}\setminus R_v$
we say that
$x$ and
$y$ are
{\bf multiplicatively-equivalent} and write $x\eqmult y$ if:
$\exists n\in\N \mbox{ s.t. }x^n \geq y \mbox{ and } y^n \geq x$.
 Note that
\begin{equation}   \label{may28}
 x\eqmult y \mbox{ if and only if } v(x) \eqad v(y)
\end{equation}
\sn
 An ordered field $K$ is an
{\bf exponential field}
if there exists a map
\[\mbox{exp}:(K,+,0,<) \longrightarrow (K^{>0},\cdot,1,<)\]
such that exp is an isomorphism of ordered groups. A map exp with
these properties will be called an {\bf exponential} on $K$. A
{\bf logarithm} on $K$ is the compositional inverse
$\mbox{log} = \mbox{exp}^{-1}$ of an exponential. Without loss of
generality, we shall always require the exponentials (logarithms)
under consideration to be
{\bf $v$-compatible}:
$\mbox{exp}(R_v)=U_v ^{>0}$ or $\mbox{log}(U_v^{>0})=R_v$.
\sn
We are mainly interested in exponentials
satisfying the {\bf growth axiom} scheme:
\mn
{\bf (GA)}\ \  $x\geq n^2\;\Longrightarrow\; \mbox{exp}(x)>x^n
\hspace{4cm} (n\geq 1)$
\mn
Note that because of the hypothesis $x\geq n^2$,
{\bf (GA)} is only relevant for $v(x) \leq 0$.
Let us consider the case $v(x)<0$. In this case,``$x>n^2\,$''
holds for all $n\in\N$ if $x$ is positive. Restricted to $K
\setminus R_v\,$, axiom scheme {\bf (GA)} is thus equivalent to the assertion
\begin{equation}                            \label{f(x)>x^n}
\forall n\in\N:\; \mbox{exp}(x)>x^n
\hspace{1cm}\mbox{ for all } x\in K^{>0}\setminus R_v\;.
\end{equation}
Applying the logarithm $\mbox{log}=\mbox{exp}^{-1}$ on both sides, we
find that this is equivalent to
\begin{equation}                            \label{increll}
\forall n\in\N:\; x> \mbox{log}(x^n)=n\mbox{log}(x)
\hspace{1cm}\mbox{ for all } x\in K^{>0}\setminus R_v\;.
\end{equation}
Via the natural valuation $v$, this in turn is equivalent to
\begin{equation}                            \label{stronglog}
v(x)<v(\mbox{log}(x)) \hspace{1cm}\mbox{ for all } x\in K^{>0}\setminus
R_v\;.
\end{equation}
A logarithm $\mbox{log}$ will be called a {\bf (GA)-logarithm}
if it satisfies (\ref{stronglog}). For more details about ordered
exponential fields and their natural valuations see [K].
\bn
In this paper, we will mainly work with ordered abelian groups and
ordered fields of the following form:
let $\Gamma$ be any totally ordered set and $R$ any ordered abelian
group. Then $R^\Gamma$ will denote the Hahn product with index
set $\Gamma$ and components $R$. Recall that this is the set of all
maps $g$ from $\Gamma$ to $R$ such that the
{\bf support}
$\{\gamma\in \Gamma\mid g(\gamma)\ne 0\}$ of $g$ is well-ordered in
$\Gamma$.
Endowed with the lexicographic order and pointwise addition, $R^\Gamma$
is an ordered abelian group, called the {\bf Hahn group}.
\sn
We want a convenient representation for the elements $g$ of the Hahn
groups. Fix a
strictly positive element $1 \in R$
(if $R$ is a field, we take $1$ to be the neutral element for
multiplication). For every $\gamma\in \Gamma$, we will denote by
$1_\gamma$ the map which sends $\gamma$ to $1$ and
every other element to $0$ ($1_\gamma$ is the characteristic function
of the singleton $\{\gamma\}$.)
 Hence, every
 $g\in R^\Gamma$ can be written in the form $\sum_{\gamma\in
\Gamma} g_{\gamma} 1_\gamma$ (where $g_{\gamma}:=g(\gamma)\in R$).
Note that
$g\eqad g'$ if and only if $\min \support g = \min \support g'$.
\sn
For $G\ne 0$ an ordered abelian group, $k$ an archimedean ordered field,
$k((G))$ will denote the (generalized) {\bf power series
field} with coefficients in $k$ and exponents in $G$. As an ordered
abelian group, this is just the Hahn group $k^G$. When we work
in $K=k((G))$, we will write $t^g$ instead of $1_g\,$. Hence, every
  series $s\in k((G))$ can be written in the form $\sum_{g\in G} s_g
t^g$ with $s_g\in k$ and well-ordered support
$\{g\in G\mid s_g\ne 0\}$.
Multiplication is given by the usual formula for multiplying series.
\sn
The natural valuation on $k((G))$ is given by
$v(s) = \min \support s$ for any series $s\in
k((G))$. Clearly the value group is (isomorphic to) $G$ and the residue
field is (isomorphic to) $k$. The valuation ring $k((G^{\geq
0}))$ consists of the
series with non-negative exponents, and the valuation ideal
$k((G^{>0}))$
of
the
series with positive exponents. The {\bf constant term} of a series $s$
is the coefficient $s_0$. The units of $k((G^{\geq 0}))$ are the series
in $k((G^{\geq 0}))$ with a non-zero constant term.
\sn
Given any series, we can
truncate it at its constant term and write it as the sum of two
series, one with strictly negative exponents, and the other with
non-negative exponents. Thus a complement in
$(k((G)), +)$ to the valuation ring is the Hahn group
$k^{G^{<0}}$. We call it the
{\bf canonical complement to the valuation ring}
and denote it by
{\bf Neg $k((G))$} or by $k((G^{<0}))$. Note that
$\mbox{ Neg }k((G))$ is in fact a (non-unital) subring , and a
$k$-algebra.
\sn
Given $s\in k((G))^{>0}$, we can
factor out the monomial of smallest exponent $g\in G$ and write
$s=t^gu$ with $u$ a unit with a positive constant term.
Thus a complement in $(k((G))^{>0}, \cdot)$ to the subgroup $U_v^{>0}$
of positive units is the group consisting of the
(monic) monomials $t^g$. We call it the
{\bf canonical complement to the positive units}
 and denote it by
{\bf Mon $k((G))$}.
\mn
Throughout this paper, {\bf fix a regular
uncountable cardinal $\kappa$}. We are particularly
interested in the
{\bf $\kappa$-bounded Hahn group} $(R^\Gamma)_{\kappa}$, the subgroup of
$R^\Gamma$ consisting of all maps of which support has cardinality
$<\kappa$.
Similarly, we consider
 the
{\bf $\kappa$-bounded power series field} $k((G))_{\kappa}$, the
subfield of
$k((G))$ consisting of all series of which support has cardinality
$<\kappa$. It is a valued subfield of $k((G))$. We denote by $k((G^{\geq
0}))_{\kappa}$ its valuation ring. A subfield
$F$ of
$k((G))$ is said to be {\bf truncation
closed} if whenever $s\in F$, then all truncations (initial segments) of
$s$ belong to
$F$ as well. If $F$ is truncation closed, then
$\mbox {Neg} (F):=\mbox { Neg } k((G))\cap F$ is a complement to
the
valuation ring of $F$.
If $F$ contains the subfield $k(t^g\
;\ g\in G)$
generated by the monic monomials, then $\mbox {Mon} (F)=\{t^g\ ;\ g\in
G\}$
is a complement to the group of positive units in $(F^{>0}, \cdot)$.
 Note that
$k((G))_{\kappa}$ is truncation
closed and contains $k(t^g\ ;\ g\in G)$.
We denote $\Neg k((G))_{\kappa}$ by
 $k((G^{<0}))_{\kappa}$.
\mn
Our goal is to define an exponential (logarithm) on $k((G))_{\kappa}$
(for appropriate choice of $G$). From the above discussion,
we get the following useful result:
\begin{proposition} \label{dec+}
 Set $K=k((G))_{\kappa}$. Then
 $(K,+,0,<)$ decomposes lexicographically as the sum:
\begin{equation}                            \label{decomp+}
(K,+,0,<) =
k((G^{<0}))_{\kappa}
\oplus k((G^{\geq 0}))_{\kappa}
\;.
\end{equation}
Similarly, $(K^{>0},\cdot,1,<)$ decomposes
lexicographically as the product:
\begin{equation}                            \label{decompx}
(K^{>0},\cdot,1,<) = \Mon (K)\times U_v^{>0}
\end{equation}
Moreover, $\Mon(K)$ is  order  isomorphic to $G$
through the isomorphism $(-v)(t^g)=-g$.
\end{proposition}
\mn
Proposition \ref{dec+} allows us to achieve our goal in two main
steps; by defining the logarithm first on $\Mon (K)$ (Lemma \ref{mainl})
and then on
$U_v^{>0}$ (Proposition \ref{an}).
%
%
\section{The Main Lemma.}          \label{main}
 We are interested in
developing a method to construct a {\bf left logarithm}
on $\R((G))_{\kappa}$, that is, an isomorphism of ordered
groups from $\Mon \R((G))_{\kappa}$ onto
$\Neg \R((G))_{\kappa}=\R((G^{<0}))_{\kappa} $.
Moreover, we want a
criterion to obtain a {\bf(GA)-left logarithm}, that is, a left
logarithm
which satisfies $t^g> \log((t^g)^n)=n\log(t^g)$ for all $n\in\N$ and
$g\in G^{<0}$.
\begin{lemma}  \label{mainl}
Let $\Gamma$ be a chain.
Set \[G:=(\R ^{\Gamma})_{\kappa} \mbox{ and } K:=\R((G))_{\kappa}\>.\]
Every isomorphism of chains \[\iota: \Gamma \ra  G^{<0}\]
lifts to an isomorphism of ordered groups
\[\hat{\iota}: (G, +) \ra  (\Neg(K), +)\]
given by
\begin{equation}  \label{hat}
\hat{\iota}(\sum_{\gamma\in\Gamma} g_{\gamma} 1_\gamma):=
\sum_{\gamma\in\Gamma} g_{\gamma} t^{\iota(\gamma)}\>
\end{equation}
 for $g=\sum_{\gamma\in\Gamma} g_{\gamma} 1_\gamma
\in G$. Furthermore, setting
\begin{equation}  \label{log}
\log (t^g)
:=\hat{\iota}(-g)=\sum_{\gamma\in\Gamma} -g_{\gamma}
t^{\iota(\gamma)}
\end{equation}
defines
 a left logarithm on $K$, which satisfies
\begin{equation} \label{may27}
v(\log t^g)=
\iota(\min\support g)
\end{equation}
 Moreover
  $\log$ is a (GA)-left
logarithm if and only if
\begin{equation}                            \label{incres}
\iota(\min\support g) > g
\hspace{1cm}\mbox{ for all } g\in G^{<0}\;.
\end{equation}
\end{lemma}
\begin{proof}
The map $\hat{\iota}$ is well defined (because of the condition imposed
simultaneously on the supports
of elements of $G$ and of $K$). It is straightforward to verify that
$\hat{\iota}$ is an isomorphism of ordered groups
and that (\ref{log}) defines a left logarithm. Also
(\ref
{incres})
follows from  (\ref {stronglog}).
\end{proof}
\begin{remark}       \label{pre}
If $\iota$ is only an embedding, one would still obtain by (\ref{hat})
an embedding
$\hat{\iota}$, and by (\ref{log}) an embedding of $\Mon(K)$ into
$\Neg(K)$ (a so called left {\it pre}-logarithm).
The maps $\hat{\iota}$ and $\log$
are surjective
(isomorphisms) if and only if
$\iota$ is surjective. This observation is used to construct
pre-logarithms on Exponential-Logarithmic Power Series fields in [K].
In this paper, we will not make use of pre-logarithms.
\end{remark}
\section{ The $\kappa$-th iterated lexicographic power of a chain.}
\label{iterate}
Let $\Gamma _0 \ne \emptyset$ be a given chain. We shall construct
canonically over $\Gamma _0$ a chain $\Gamma _{\kappa}$ together with an
isomorphism of ordered chains
\[\iota_{\kappa}: \Gamma _{\kappa}\ra
G_{\kappa}^{<0}\] where
$G_{\kappa}:=(\R^{\Gamma_{\kappa}})_{\kappa}$.
We call the pair $(\Gamma _{\kappa}, \iota_{\kappa})$ the $\kappa$-th
{\bf iterated lexicographic power} of  $\Gamma _0$.
\sn
We shall construct by transfinite induction on $\mu \leq \kappa$ a chain
$\Gamma _{\mu}$
together with an
embedding of ordered chains
\[\iota_{\mu}: \Gamma _{\mu}\ra
G_{\mu}^{<0}\] where
$G_{\mu}:=(\R^{\Gamma_{\mu}})_{\kappa}$. We shall have
$\Gamma _{\nu}\subset \Gamma _{\mu}$ and
$ \iota_{\nu}\subset \iota _{\mu}$ if $\nu < \mu$.
\sn
For $\mu=0$, set
$G_0=(\R^{\Gamma_0})_{\kappa}$ and
$\iota _0\;:\;\Gamma_0 \ra  G_0^{<0}$ be defined by $\gamma \mapsto
-1_{\gamma}$.
Now assume that for all $\alpha < \mu$ we have already constructed
$\Gamma_{\alpha}$,
$G_{\alpha}:=(\R^{\Gamma_{\alpha}})_{\kappa}$,
 and the embedding
\[\iota_{\alpha}: \Gamma _{\alpha}\ra G_{\alpha}^{<0}\>.\]
First assume that $\mu=\alpha + 1$ is a successor ordinal.
Since $\Gamma _{\alpha}$ is isomorphic to a subchain of $G_{\alpha}^{<0}$
through $\iota_{\alpha}$, we can take $\Gamma _{\alpha+1}$ to be a chain
containing $\Gamma _{\alpha}$
as a subchain and admitting an isomorphism $\iota_{\alpha+1}$ onto
$G_{\alpha}^{<0}$ which extends
$\iota_{\alpha}\,$. More precisely,
\[\Gamma _{\alpha+1}:=\Gamma _{\alpha}\cup
(G_{\alpha}^{<0}\setminus \iota_{\alpha}(\Gamma _{\alpha}))\>,\]
endowed with the {\bf patch ordering}: if
$\gamma_1$, $\gamma_2 \in \Gamma _{\alpha+1}$ both belong to
$\Gamma_{\alpha}$, compare them there, similarly if they both belong to
$G_{\alpha}^{<0}$. If $\gamma_1 \in \Gamma_{\alpha}$ but
$\gamma_2 \in G_{\alpha}^{<0}$ we set $\gamma_1 < \gamma_2$ if and only
if $\iota_{\alpha}(\gamma_1) < \gamma_2$ in $G_{\alpha}$.
Then $\iota_{\alpha+1}$ is defined in the obvious way:
$\iota_{\alpha+1}|_{\Gamma _{\alpha}}:= \iota_{\alpha}$ and
$\iota_{\alpha+1}
|_{(G_{\alpha}^{<0}\setminus \iota_{\alpha}(\Gamma _{\alpha}))}:=$ the
identity map. Note that
\begin{equation}  \label{onto}
\iota_{\alpha+1}(\Gamma _{\alpha+1})= G_{\alpha}^{<0}.
\end{equation}
Thus $\iota_{\alpha+1}$ is an embedding of $\Gamma _{\alpha+1}$
into $G_{\alpha+1}^{<0}$.
\n
If $\mu$ is a limit ordinal we set
\[\Gamma _{\mu}:=\bigcup_{\alpha <
\mu} \Gamma_{\alpha}\>\>,\>\>
\iota_{\mu}:=
\bigcup_{\alpha <\mu} \iota_{\alpha}\;
\mbox{ \ \ and \ \ }
G_{\mu}:=(\R^{\Gamma_{\mu}})_{\kappa}\>.\]
Note that
by construction and (\ref{onto})
\begin{equation} \label{ontolim}
\iota_{\mu}(\Gamma _{\mu})= \bigcup_{\alpha <\mu}G_{\alpha}^{<0}
\end{equation}
 and
$\bigcup_{\alpha <\mu} G_{\alpha}\subset G_{\mu}$.
\sn
This completes the construction of
$\Gamma _{\kappa}:=\bigcup_{\alpha <\kappa} \Gamma_{\alpha}\>,$
$\iota_{\kappa}:=
\bigcup_{\alpha <\kappa} \iota_{\alpha}$
and
$G_{\kappa}:=(\R^{\Gamma_{\kappa}})_{\kappa}$. We now claim that
\[G_{\kappa}=\bigcup_{\alpha <\kappa}G_{\alpha}\]
(Once the claim is established, we conclude from (\ref{ontolim}) that
$\iota_{\kappa}:\Gamma _{\kappa}\ra G_{\kappa}^{<0}$ is an isomorphism,
as required). Let $g\in G_{\kappa}$ and $\kappa>\delta:=\card(\support
g)$. Now $\support g:=\{\gamma_{\mu}\>;\> \mu<\delta\}\subset
\Gamma_{\kappa}$, so for
every $\mu < \delta$ choose $\alpha_{\mu} < \kappa$ such that
$\gamma_{\mu}\in \Gamma_{\alpha_{\mu}}$. Clearly
$\card(\{\alpha_{\mu}\>;\> \mu<\delta\})\leq \delta < \kappa$ so
$\{\alpha_{\mu}\>;\> \mu<\delta\}$ cannot be cofinal in $\kappa$ (since
$\kappa$ is regular), therefore it is bounded above by some
$\alpha\in \kappa$. It follows that $\support g \subset
\Gamma_{\alpha}$, so $g\in G_{\alpha}$ as required.
\begin{proposition}    \label{twist}
Assume that $\sigma\in \Aut(\Gamma_{\kappa})$  is such that
$\sigma |_{ \Gamma_{\mu}}\in \Aut(\Gamma_{\mu})$ for all $\mu\in
\kappa$ and
$\sigma(\gamma) > \gamma$ for all $\gamma\in \Gamma_0$.
Then the isomorphism
\[l:=\iota_{\kappa}\circ \sigma: \Gamma _{\kappa}\ra G_{\kappa}^{<0}\>\]
satisfies (\ref{incres}).
\end{proposition}
\begin{proof}
Let $g\in G_{\kappa}^{<0}$ and $\gamma_{\mu}:=\min\support g\in
\Gamma_{\mu}$ for the least such $\mu\in \kappa$. We
prove that (\ref{incres}) holds by transfinite induction on $\mu$.
If $\mu=0$, then $\gamma_0 \in \Gamma_0$ so
\[l(\gamma_0)=\iota_0\circ \sigma (\gamma_0)=-1_{\sigma(\gamma_0)}
> g\>.\]
Now assume that the assertion holds for all $\alpha < \mu$. Since
\[\iota_{\kappa}\circ \sigma (\Gamma_{\alpha+1})=
\iota_{\alpha+1}(\Gamma_{\alpha+1})= G_{\alpha}^{<0}\>,\]
by (\ref {onto}) and for $\mu$ limit
\[\iota_{\kappa}\circ \sigma (\Gamma_{\mu})=
\iota_{\mu}(\Gamma _{\mu})= \bigcup_{\alpha <\mu}G_{\alpha}^{<0}\]
by (\ref{ontolim}),
we have in any case that
\begin{equation}   \label{ontol}
l(\gamma_{\mu})\in  G_{\alpha}^{<0} \mbox{ for some }
\alpha<\mu .
\end{equation}
Set $l(\gamma_{\mu}):=g'\in  G_{\alpha}^{<0}$.
 We have to show that $g<g'$, for this it is enough to show
that
$\min\support g<
\min\support g'$, or equivalently that:
\[l(\min\support g)<l(\min\support g')\>.\]
But the last inequality holds since by induction assumption we have that
$g' < l(\min\support g')$.
\end{proof}
\begin{proposition} \label {trevor}
Let $\sigma _0\in \Aut(\Gamma_0)$.
Then $\sigma _0$ can be extended to $\sigma \in
\Aut(\Gamma_{\kappa})$ satisfying $\sigma |_{\Gamma_{\mu}}\in
\Aut(\Gamma_{\mu})$ for all $\mu\in \kappa$.
In particular, if $\sigma _0\in \Aut(\Gamma_0)$
satisfies $\sigma_0(\gamma)>\gamma$ for all $\gamma\in
\Gamma_0$, then $\sigma$ satisfies the hypothesis of Proposition
\ref{twist}.
\end{proposition}
\begin{proof}
We first note that any $\sigma _{\mu}\in
\Aut(\Gamma_{\mu})$ lifts to
$\hat{\sigma} _{\mu}\in \Aut(G_{\mu})$ as follows.
For $g=\sum_{\gamma\in\Gamma_{\mu}} g_{\gamma} 1_\gamma
\in G_{\mu}$, set:
\begin{equation}  \label{hataut}
\hat{\sigma} _{\mu}(\sum_{\gamma\in\Gamma_{\mu}} g_{\gamma} 1_\gamma):=
\sum_{\gamma\in\Gamma_{\mu}} g_{\gamma} 1_{\sigma _{\mu}(\gamma)}\>
\end{equation}
Observe that if $\alpha < \mu$ and
$\sigma _{\mu}\in
\Aut(\Gamma_{\mu})$ extends
$\sigma _{\alpha}\in
\Aut(\Gamma_{\alpha})$, then also $\hat{\sigma} _{\mu}$ extends
$\hat{\sigma}_{\alpha}$.
By induction on $\mu \leq \kappa$, we now construct
$\sigma _{\mu}\in
\Aut(\Gamma_{\mu})$ satisfying the following two properties:
\begin{equation}  \label{commuting}
(i)\ \ \hat{\sigma} _{\mu}\circ \iota_{\mu} = \iota_{\mu} \circ
\sigma_{\mu}
\
\ \mbox { and } \ \
(ii)\ \ \sigma _{\mu} \supset
\sigma _{\beta}\ \  \mbox{ for all } \beta \leq \mu.
\end{equation}
Note that (15) part (i) implies that
\begin{equation} \label{welldefined}
\mbox{ for all } g\in
G_{\mu}^{<0}: \ \ \ \
 \hat{\sigma} _{\mu}(g) \in \iota_{\mu}(\Gamma_{\mu})\ \  \mbox{ if and
only if }\ \ g\in \iota_{\mu}(\Gamma_{\mu})
\end{equation}
 It is readily verified that $\sigma _0$
satisfies
(\ref{commuting}).
Assume that for $\alpha < \mu$, $\sigma _{\alpha}$ has been
constructed
satisfying
(\ref{commuting}).\n
If $\mu=\alpha +1$, define $\sigma
_{\alpha+1}$ on
$\Gamma _{\alpha+1}=\Gamma _{\alpha}\cup
(G_{\alpha}^{<0}\setminus \iota_{\alpha}(\Gamma _{\alpha}))\>$ by
setting:
$\sigma_{\alpha+1}|_{\Gamma _{\alpha}}:= \sigma_{\alpha}$ and
$\sigma_{\alpha+1}
|_{(G_{\alpha}^{<0}\setminus \iota_{\alpha}(\Gamma _{\alpha}))}
:=\hat{\sigma}_{\alpha}$. Since $\hat{\sigma} _{\alpha}$
satisfies (\ref{welldefined}),
$\sigma _{\alpha+1}$ is well--defined. It easily follows from the
definition of $\sigma _{\alpha+1}$
that
$\sigma _{\alpha+1}\supset \sigma _{\alpha}$,
and that
$\sigma _{\alpha+1}$ is a bijection satisfying (\ref{commuting}).
It remains to verify that $\sigma _{\alpha+1}(\gamma_1)
< \sigma _{\alpha+1}(\gamma_2)$ for $\gamma_1 < \gamma_2$,
$\gamma_1, \gamma_2 \in \Gamma_{\alpha +1}$. We only verify this when
$\gamma_1\in
\Gamma
_{\alpha}$ and
$\gamma_2\in G_{\alpha}^{<0}$ (the verification in the other cases is
straightforward).
\relax From $\iota_{\alpha}(\gamma_1) < \gamma_2$ in $G_{\alpha}$ follows that
$\hat{\sigma}_{\alpha}(\iota_{\alpha}(\gamma_1)) <
\hat{\sigma}_{\alpha}(\gamma_2)$ in $G_{\alpha}$. By (\ref{commuting}),
we therefore have
$\iota_{\alpha}(\sigma_{\alpha}(\gamma_1)) <
\hat{\sigma}_{\alpha}(\gamma_2)$ in $G_{\alpha}$. That is,
$\iota_{\alpha}(\sigma_{\alpha+1}(\gamma_1)) <
\sigma_{\alpha+1}(\gamma_2)$ in $G_{\alpha}$, or equivalently
$\sigma_{\alpha+1}(\gamma_1)<
\sigma_{\alpha+1}(\gamma_2)$ in $\Gamma_{\alpha+1}$ as required.\n
Finally, if $\mu$ is a limit ordinal, set $\sigma_{\mu}
:=
\bigcup_{\alpha <\mu} \sigma_{\alpha}$.
Then
$\sigma:=\sigma _{\kappa}$ is the required
$\sigma \in
\Aut(\Gamma_{\kappa})$.
\end{proof}
%
%
%
\section{$\kappa$--bounded models.}          \label{models}
We now extend the definition of the logarithm to the positive units. Below,
for $r\in\R$, $r>0$ we denote by $\log r$ the natural logarithm of $r$.
\begin{proposition}  \label{an}
Let $G$ be any divisible ordered abelian group, and set
$K:=\R((G))_{\kappa}$. For $u\in U_v^{>0}$ write $u= r(1+\varepsilon)$
(with $r\in\R$, $r>0$ and $\varepsilon \in I_v$ infinitesimal).
Then
\begin{equation}                            \label{nl}
\log(u):=\log r(1+\varepsilon)  = \log r + \sum_{i=1}^{\infty}
(-1)^{(i-1)}
\frac{\varepsilon^i}{i}
\end{equation}
defines an isomorphism of ordered groups from $U_v^{>0}$ onto $R_v$
\end{proposition}
\begin{proof}
The formal sum given in (\ref{nl}), and more generally, any formal sum
$\sum_{i=0}^{\infty}
r_i \varepsilon^i$ (with $r_i\in\R$)
  is
a well-defined element of
$\R((G))$:
it has well-ordered support, since $\support
\varepsilon \subset G^{>0}$. Also, the map defined by (\ref{nl}) is a
bijective, order preserving group homomorphism
cf.\ [F]. It remains to verify that
\[\card(\support \varepsilon) < \kappa\Longrightarrow \card(\support
\sum_{i=0}^{\infty} r_i \varepsilon^i) < \kappa \>.\]
Note that \[\support r_i \varepsilon^i\subset \oplus_i
\support \varepsilon:=\{g_1+\cdots + g_i\>|\> g_j \in \support
\varepsilon \mbox { for all } j=1,\cdots,i\}\>,\]
 and clearly, $\card(\oplus_i \support
\varepsilon)< \kappa$ for all $i$,
so $\card(\cup_i (\oplus_i \support \varepsilon)) <
\kappa$.
Now observe that
$\support \sum_{i=0}^{\infty} r_i \varepsilon^i\subset \cup_i (\oplus_i
\support \varepsilon)\>.$
\end{proof}
We can now define the logarithm on the
positive elements of $\R((G_{\kappa}))_{\kappa}$
making $\R((G_{\kappa}))_{\kappa}$ into a model of $T_{\exp}$:= {\bf the
elementary theory of the reals with exponentiation}. Below,
$T_{\rm an}$ := the theory of the reals with restricted
analytic functions and
$T_{{\rm an} , \exp}$ := the theory of the reals with restricted
analytic functions and exponentiation (see [D--M--M1] for
axiomatizations of these theories).
\begin{theorem}        \label{T(an,exp)}
 Let $\kappa$ be a regular uncountable cardinal, $\Gamma_0$ a chain,
$\Gamma_{\kappa}$
the $\kappa$-th
lexicographic iterated power of $\Gamma_0$, and
$G_{\kappa}=(\R^{\Gamma_{\kappa}})_{\kappa}$.
Let $\sigma\in \Aut(\Gamma_{\kappa})$
 and
\[l: \Gamma _{\kappa}\ra G_{\kappa}^{<0}\>\]
be as in Proposition \ref{twist}.
For positive $a \in \R((G_{\kappa}))_{\kappa}$, write
$a=t^gr(1+\varepsilon)$, with
$g=\sum_{\gamma\in\Gamma_{\kappa}} g_{\gamma} 1_\gamma
\in G_{\kappa}$, $r\in \R^{>0}$, and $\varepsilon$ infinitesimal.
Then
\begin{equation}  \label{logall}
\log (a):=\log (t^gr(1+\varepsilon))=\sum_{\gamma\in\Gamma} -g_{\gamma}
t^{l(\gamma)} +
\log r + \sum_{i=1}^{\infty}
(-1)^{(i-1)} \frac{\varepsilon^i}{i}
\end{equation}
defines a logarithm on $\R((G_{\kappa}))_{\kappa}^{>0}$
making $\R((G_{\kappa}))_{\kappa}$ into a model of $T_{\exp}$.
\end{theorem}
\begin{proof}
By Lemma \ref{mainl}, Proposition \ref{twist}, and Proposition
\ref{an}, (\ref{logall}) defines a (GA)-logarithm.
Using the Taylor expansion of any analytic function, one can endow
$\R((G_{\kappa}))_{\kappa}$ with a natural interpretation
of the restricted analytic functions
(as we did in
Proposition \ref{an} for the logarithm). This makes
 $\R((G_{\kappa}))_{\kappa}$
into a substructure of the $T_{\rm an}$ model $\R((G_{\kappa}))$
(cf. [D--M--M1]).
\relax From the quantifier elimination results of [D--M--M1], we get that
$\R((G))_{\kappa}$ is a model of $T_{\rm an}$. Since log is a (GA)-logarithm,
it follows (from the axiomatization given in [D--M--M1])
that $\R((G))_{\kappa}$ is a model of $T_{{\rm an},\exp}$ .
\end{proof}
%
%
%
\section{Growth Rates.}          \label{rank}
Let $\Gamma$ be a chain and $\sigma\in \Aut(\Gamma)$. Assume that
\begin{equation}    \label{incaut}
\sigma(\gamma)> \gamma \mbox{   for all   } \gamma\in \Gamma
\end{equation}
An automorphism satisfying (\ref{incaut}) will be called an {increasing
automorphism}. By induction, we define the {\bf n-th iterate}
of $\sigma$: $\sigma^1(\gamma):=\sigma
(\gamma)$ and
$\sigma^{n+1}(\gamma):=\sigma(\sigma^n(\gamma))$.
We define an equivalence relation on $\Gamma$ as follows:
For $\gamma,\gamma'\in\Gamma$, set
\begin{equation} \label{sigmaeq}
\gamma\sim_{\sigma}\gamma' \mbox{  if and only  } \exists n\in\N \mbox {
such that }
\sigma^n (\gamma)\geq \gamma'  \> \mbox{ and } \>
\sigma^n (\gamma')\geq \gamma\
\end{equation}
The
equivalence classes $[\gamma]_{\sigma}$ of $\sim_{\sigma}$ are
convex and closed under application of $\sigma$. By the convexity, the
order of $\Gamma$ induces an order on $\Gamma/{\sim}_{\sigma}$ such
that
$[\gamma]_{\sigma}<[\gamma']_{\sigma}$ if $\gamma <
\gamma'$. The order type of $\Gamma/{\sim}_{\sigma}$ is the {\bf rank}
of $(\Gamma, \sigma)$.
\mn
Similarly, let $K$ be a real closed field and log a ({\bf GA})-
logarithm on
$K^{>0}$. Define an equivalence relation on  $K^{>0}\setminus R_v$:
\begin{equation}
a\sim_{log}a' \mbox{ if and only if }
\exists n\in\N \mbox{ such that }
 \log_n(a)\leq (a') \> \mbox{ and } \>  \log_n(a')\leq a
\end{equation}
(where $\log_{n}$ is the n-th iterate of the log).
Again, the
log-equivalence classes are
convex and closed under application of log.
 The order type of the
chain of equivalence classes is the
{\bf
logarithmic rank} of $(K^{>0}, \log)$.
Note that if
 $x$ and $y$ are archimedean-equivalent or multiplicatively-equivalent
(cf.\ (\ref{may28})),
then they are {\it a fortiori} log-equivalent.
\mn
We now compute the logarithmic rank of the
models described in Theorem \ref{T(an,exp)}. Below, set
$\sigma _0:= \sigma |_{ \Gamma _0}$.
\begin{theorem}      \label{preservation}
The logarithmic rank of $(\R((G_{\kappa}))_{\kappa}^{>0}, \log)$ is
equal to the rank of $(\Gamma_0, \sigma _0)$.
\end{theorem}
\begin{proof}
Let $a \in K^{>0}\setminus R_v$ , write $a=t^gu$ (with $u$ a unit, $g\in
G_{\kappa}^{<0}$).
Since $a$ is archimedean-equivalent to $t^g$, it is log-equivalent to
it. So it is enough to consider monomials  $t^g$ with
 $g=\sum_{\gamma\in\Gamma_{\kappa}} g_{\gamma} 1_\gamma \in
G_{\kappa}^{<0}$.
 Set
$\gamma_{\mu}:=\min \support g \in \Gamma_{\mu}$ for the least such
$\mu\in
\kappa$.
We show by transfinite induction on $\mu$ that there exists $g_0\in
G_{\kappa}^{<0}$ such that $\gamma_0:=\min \support g_0 \in
\Gamma_0$ and
$t^g$ is log-equivalent to $t^{g_0}$.
\sn
If $\mu=0$ there is nothing to prove. Assume that the assertion holds
for all $\alpha<\mu$. Now
\begin{equation}  \label{logl}
\log (t^g)=\sum_{\gamma\in\Gamma} -g_{\gamma}
t^{l(\gamma)}
\end{equation}
 is archimedean-equivalent (cf.\ (\ref{may27})),
so log-equivalent to $t^{l(\gamma_{\mu})}$. By (\ref{ontol}) and
induction hypothesis, the assertion holds for $t^{l(\gamma_{\mu})}$,
and thus for $t^{g}$  by transitivity.
\sn
Now we determine the logarithmic equivalence class
of $t^g$ for $g\in G_{\kappa}^{<0}$
such that
$\gamma_0:=\min \support g \in
\Gamma_0$. Now $t^g$ is multiplicatively-equivalent,
so log-equivalent to $t^{-1_{\gamma_0}}$, so it is enough to consider
monomials of the form $t^{-1_{\gamma}}$ with $\gamma\in \Gamma_0$.
We claim that
\[\mbox{ for all } \gamma, \gamma'\in \Gamma_0:
t^{-1_{\gamma}} \sim_{log} t^{-1_{\gamma'}}
\mbox{ if and only if }\gamma\sim_{\sigma}\gamma'\>.\]
We first find a formula for $\log_n(t^{-1_{\gamma}} )$. Using
(\ref {logl})
 we compute:
$\log (t^{-1_{\gamma}})=
t^{l(\gamma)}=t^{\iota_0 \circ \sigma(\gamma)}=
t^{\iota_0 (\sigma(\gamma))}= t^{-1_{\sigma(\gamma)}}$
(since $\sigma(\gamma) \in \Gamma_0$). By induction, we see that
for all $n\in \N$:
\[\log_n (t^{-1_{\gamma}})= t^{-1_{\sigma^n(\gamma)}}\>.\]
We conclude: $\gamma\sim_{\sigma}\gamma' \Llra \exists n\in \N$ such
that
$\sigma^n (\gamma)\geq \gamma'  \> \mbox{ and } \>
\sigma^n (\gamma')\geq \gamma \Llra 1_{\sigma^n (\gamma)} \leq
1_{\gamma'}$ and
$1_{\sigma^n (\gamma')} \leq
1_{\gamma}\Llra -1_{\gamma'}\leq -1_{\sigma^n (\gamma)}$ and
$-1_{\gamma}\leq -1_{\sigma^n (\gamma')}\Llra$
\[t^{-1_{\gamma'}} \geq
t^{-1_{\sigma^n (\gamma)}}=\log_n (t^{-1_{\gamma}})\mbox { and }
t^{-1_{\gamma}} \geq
t^{-1_{\sigma^n (\gamma')}}=\log_n (t^{-1_{\gamma'}})\>,\]
if and only if $t^{-1_{\gamma}} \sim_{log} t^{-1_{\gamma'}}$ as
required.
\end{proof}
\begin{theorem} \label{final}
Let $\kappa$ be a regular uncountable cardinal with
$\kappa=\kappa^{<\kappa}$.
Let $\Gamma_0$ be any chain of cardinality $\kappa$ which admits a family ${\cal A}=\{\sigma_0^{\alpha}\>|\> \alpha \in 2^{\kappa}\}\subset \Aut (\Gamma_0)$ of increasing automorphisms of pairwise distinct ranks. Let $\Gamma_{\kappa}$ be
the $\kappa$-th iterated lexicographic power of  $\Gamma _0$,
$G_{\kappa}:=(\R^{\Gamma_{\kappa}})_{\kappa}$ the corresponding
$\kappa$-bounded Hahn group, and $K=\R((G_{\kappa}))_{\kappa}$
the corresponding $\kappa$-bounded power series field of cardinality
$\kappa$. Then $K$ admits a family $\{\exp ^{\alpha}\>|\> \alpha\in
2^{\kappa}\}$ of
 $2^{\kappa}$ exponentials.
For every $\alpha\in 2^{\kappa}$,
$(K,\exp^ {\alpha})$ is a model of real exponentiation.
The $2^{\kappa}$ exponentials are of pairwise distinct exponential rank,
but all agree on the valuation ring of $K$.
\end{theorem}
\begin{proof}
For every $\sigma_0^{\alpha}$, let
$\sigma ^ {(\alpha)} \in \Aut (\Gamma_{\kappa})$ be the corresponding
extension
(Proposition \ref {trevor}). Set
$l^{\alpha}:=\iota_{\kappa}\circ \sigma ^ {(\alpha)}$,
and let $\log ^{\alpha}$
be the corresponding logarithm (obtained by replacing in $l$
by $l^{\alpha}$ in equation (\ref{logall}) ).
Now apply Theorem \ref{preservation}.
\end{proof}
In the next section, we give an explicit construction
of chains satisfying the hypothesis of this theorem.
%
%
%
\section{Chains with $2^{\kappa}$ automorphisms of
distinct ranks.}  \label{marker}

\begin{lemma}                                      \label{tor1}
Let $\beta$ be an ordinal, and consider the chain $\Gamma_0:= \beta
\lex\Q$ . For every $\alpha\in\beta$, let $\Q_{\alpha}$,
be the $\alpha$th-copy of $\Q$. Fix
$\tau_{\alpha}$ and $\tau'_{\alpha} \in\Aut(\Q_{\alpha})$
increasing automorphisms of rank $1$ and $\Z$ respectively. For every
$S\subset \beta$ define $\tau_S$ as follows:
\[\tau_S|_{\Q_{\alpha}}   := \left\{\begin{array}{ll}
 \tau_{\alpha}& \mbox{if } \alpha\in S \\
\tau'_{\alpha} & \mbox{otherwise.}
\end{array}\right.\]
Then the rank of $\tau_S=\sum _{\alpha \in \beta} \delta_S(\alpha)$,
where
\[\delta_S(\alpha)
 := \left\{\begin{array}{ll}
 1& \mbox{if } \alpha\in S \\
\Z & \mbox{otherwise.}
\end{array}\right.\]
\end{lemma}
Lemma $\ref{tor1}$ is a consequence of the following more general
observation:
\begin{proposition}   \label{tor2}
Let $I$ be a chain, and $\{(\Gamma_i , \tau_i) \> | \> i\in I\}$ a collection of chains
$\Gamma_i$ endowed with an increasing automorphism $\tau_i$. Set
\[\Gamma :=\sum_{i\in I} \Gamma_i \mbox{ and } \tau :=\sum_{i\in I}
\tau_i\>,\]
(that is, $\tau|_{\Gamma_i}=\tau_i$).
Then the rank of $(\Gamma, \tau)$ is equal to $\sum_{i\in I} \mbox{ rank
} (\Gamma_i , \tau_i)$.
\end{proposition}
The proof is straightforward and we omit it.
\begin{remark}
(i)\ \ In [H--K--M], other
arithmetic operations on chains are studied; it may be interesting
for future work, to study the behaviour of automorphism ranks with
respect to these operations.
\sn
(ii)
Automorphisms $\tau_{\alpha}$ and $\tau'_{\alpha} \in\Aut(\Q_{\alpha})$
such as in Lemma \ref {tor1} exist: for example, set $\tau(q):=q+1$,
$\tau\in\Aut(\Q)$ is of rank
1. To produce $\tau' \in\Aut(\Q)$ of rank $\Z$, note that by Cantor's
Theorem $\Q\simeq\Z\lex \Q$. Define $\tau'$ piecewise as follows: for
$z\in\Z$ we let $\tau'|_{\Q_z}\in\Aut(\Q_z)$ be the translation
automorphism
$\tau'(q)=q+1$ for $q\in \Q_z$, then $\tau'$ is defined by patching, and
has clearly rank $\Z$ as required.
\sn
(iii)\ \ If $\beta$ is an infinite cardinal, then
$\card(\beta\lex \Q)=\beta$.
\end{remark}
We now state and prove the main result of this section. Below, we keep
the notation of Lemma \ref{tor1}.
\begin{proposition} \label{tor3}
Let $\beta$ be an ordinal and $s\subset \beta$. Set
\[\Delta_S:=\sum_{\alpha\in\beta} \delta_S(\alpha)\> .\]
Then
\[\Delta_S\simeq\Delta_{S'}\mbox{ if and only if } S=S'\>.\]
\end{proposition}
\begin{proof}
Fix an isomorphism $\varphi: \Delta_S\simeq\Delta_{S'}$. We show by
induction on $\alpha\in\beta$ that
\begin{equation} \label{tor4}
\varphi(\delta_S(\alpha)) = \delta_{S'}(\alpha) .
\end{equation}
(The Proposition is proved once (\ref{tor4}) is established:
it follows from (\ref{tor4}) that $\delta_S(\alpha))=1$ if and only if
$\delta_S'(\alpha)=1$ i.\ e.\ $ S=S'$.)
Let $\alpha=0$. Assume that $\delta_S(0)=1$. Then necessarily
$\delta_{S'}(0)=1$ and (\ref{tor4}) holds (since $\varphi$ has to map
the least element of $\Delta_S$ to the  least element of $\Delta_{S'}$).
 Assume now that $\delta_S(0)=\Z$, then necessarily $\delta_{S'}(0)=\Z$.
We claim that (\ref{tor4}) holds in this case too. Clearly, since
$\delta_S(0)$ is an initial segment of  $\Delta_S$,
$\varphi(\delta_S(0))$ is an initial segment of  $\Delta_{S'}$. It thus
suffices to show that $\varphi(\delta_S(0))\subset \delta_{S'}(0)$.
Assume for a contradiction that
$\varphi(\delta_S(0))\cap \delta_{S'}(1)\not=\emptyset$. There are 2
cases to consider. If $\delta_{S'}(1)=1$, then $1$ has left character
$\aleph_0$. This is impossible since no such element exists in
$\delta_S(0)$. If  $\delta_{S'}(1)=\Z$, then $\varphi(\delta_S(0))$ has
an $\aleph_0 \aleph_0$-gap. This is impossible since no such gap exists
in $\Z$. The claim is established.\sn Now assume that (\ref{tor4}) holds
for all $\alpha < \mu < \beta$, we show it holds for $\mu$. From
induction hypothesis we deduce that
\begin{equation} \label{tor5}
\varphi(\sum_{\alpha < \mu }\delta_S(\alpha))=\sum_{\alpha < \mu
}\delta_{S'}(\alpha) \>,
\end{equation}
therefore
\begin{equation}  \label{tor6}
\varphi(\sum_{\nu \geq \mu }\delta_S(\nu))=
\sum_{\nu \geq \mu }\delta_{S'}(\nu) \>.
\end{equation}
With the help of (\ref{tor5}) and (\ref{tor6}), the same argument as the
one used for the
induction begin  (with $\mu$ and $\mu+1$ instead of $0$ and $1$)
applies now to establish (\ref{tor4}) for $\mu$.
\end{proof}
\begin{corollary} \label{tor7}
The chain $\Gamma_0=\kappa\lex \Q$ admits of family of $2^{\kappa}$
increasing automorphisms, of pairwise distinct ranks.
\end{corollary}
\bn
\bn
\bn
{\bf References}
{\small\rm
\newenvironment{reference}%
{\begin{list}{}{\setlength{\labelwidth}{5em}\setlength{\labelsep}{0em}%
\setlength{\leftmargin}{5em}\setlength{\itemsep}{-1pt}%
\setlength{\baselineskip}{3pt}}}%
{\end{list}}
\newcommand{\lit}[1]{\item[{#1}\hfill]}
\begin{reference}
\lit{[A--K]} Alling, N.L.\ -- Kuhlmann, S.$\,$: {\it On $\eta_\alpha$-groups and fields}, Order {\bf 11} (1994), 85-92
\lit{[D--M--M1]} van den Dries, L.\ -- Macintyre, A.\ --
Marker, D.$\,$: {\it The elementary theory of restricted analytic
functions with exponentiation}, Annals Math.\ {\bf 140} (1994), 183--205
\lit{[D--M--M2]} van den Dries, L.\ -- Macintyre, A.\ -- Marker, D.$\,$:
{\it Logarithmic-Exponential series}
Annals Pure and Aplied Logic {\bf 111}
(2001), 61--113
\lit{[F]} Fuchs, L.$\,$: {\it Partially ordered algebraic
systems}, Pergamon Press, Oxford (1963)
\lit{[H--K--M]} Holland, W. C.\ -- Kuhlmann, S.\ -- McCleary, S.$\,$:
{\it Lexicographic Exponentiation of chains}, to appear in the Journal
of Symbolic Logic
\lit{[K]} Kuhlmann, S.$\,$: {\it Ordered
Exponential
Fields},
The Fields Institute Monograph Series, vol.\ 12, AMS Publications
(2000)
\lit{[K--K--S]} Kuhlmann, F.-V.\ -- Kuhlmann, S.\ -- Shelah, S.$\,$:
{\it Exponentiation in power series fields},
 Proc.\ Amer.\ Math.\ Soc.\ {\bf 125} (1997), 3177-3183
\lit{[K--T]} Kuhlmann, S.\ -- Tressl,  M.$\,$:
{\it A Note on Logarithmic - Exponential and Exponential - Logarithmic
Power Series Fields}, work in progress (2004)
\lit{[T]} Tarski, A.$\,$: {\it A Decision Method for Elementary Algebra and Geometry}, 2nd Edition, University of California Press, Berkeley, Los Angeles, CA (1951)
\lit{[W]} Wilkie, A.$\,$: {\it Model completeness results for expansions of the ordered field of real numbers by restricted Pfaffian functions and the exponential function},\n J. Amer. Math. Soc. {\bf 9} (1996), 1051--1094
\end{reference}}
\bn
\bn
\parbox[t]{8cm}{\small\rm
Research Unit Algebra and Logic\n
University of Saskatchewan\n
Mc Lean Hall, 106 Wiggins Road\n
Saskatoon, SK S7N 5E6\n
email: skuhlman@math.usask.ca}
\hfil
\parbox[t]{6.5cm}{\small\rm
Department of Mathematics\n
The Hebrew University of Jerusalem\n
Jerusalem, Israel\n
email: shelah@math.huji.ac.il}
\end{document}